\documentclass[12pt]{article}
\usepackage{amsmath,amssymb,amsthm,amscd,graphics
}
\newcommand\germa{{\mathfrak a}}
\newcommand\Ad{\operatorname{Ad}}
\newcommand\ad{\operatorname{ad}}

\newcommand\borel{{\mathfrak b}}
\newcommand\C{{\mathbb C}}
\newcommand\caff{{{\mathrm c}_{\mathrm a}}}
\newcommand\calO{{\mathcal O}}
\newcommand\CB{\operatorname{CB}}
\newcommand\CC{\operatorname{CC}}
\newcommand\eps{\varepsilon}
\newcommand\F{{\mathcal F}}
\newcommand\h{{\mathfrak h}}
\newcommand\Hom{\operatorname{Hom}\nolimits}
\newcommand\g{{\mathfrak g}}
\newcommand\gout{\g_{\mathrm{out}}}
\newcommand\gtw{\g^{\mathrm{tw}}}
\newcommand\Gtw{G^{\mathrm{tw}}}
\newcommand\Ind{\operatorname{Ind}\nolimits}
\newcommand\K{{\mathcal K}}
\newcommand\khat{{\hat k}} 
\newcommand\m{{\mathfrak m}}
\newcommand\n{{\mathfrak n}}
\newcommand\orb{\mathrm{orb}}
\newcommand\out{\mathrm{out}}
\newcommand\Proj{{\mathbb P}}
\newcommand\Res{\mathop{\operatorname{Res}}\nolimits}
\renewcommand\setminus{\smallsetminus}
\newcommand\simeqq{\cong}
\newcommand\tensor{\otimes}
\newcommand\tr{\operatorname{tr}}
\newcommand\trig{\mathrm{trig}}
\newcommand\univac{\boldsymbol{1}}
\newcommand\wt{\operatorname{wt}}
\newcommand\Z{{\mathbb Z}}
\newtheorem{theorem}{Theorem}[section]
\newtheorem{lemma}[theorem]{Lemma}

\theoremstyle{definition}
\newtheorem{definition}[theorem]{Definition}
\theoremstyle{remark}

\newcommand\thmref[1]{Theorem~\ref{#1}}

\newcommand\secref[1]{\S\ref{#1}}

\newcommand\lemref[1]{Lemma~\ref{#1}}

\newcommand\defref[1]{Definition~\ref{#1}}

\newcommand\Refcite[1]{Ref.~\cite{#1}}
\begin{document}

\markboth{Takashi TAKEBE}
{Trigonometric degeneration and orbifold WZW model. I}

%
%

\title{Trigonometric Degeneration\\
and\\
Orbifold Wess-Zumino-Witten Model. I
}

\author{Takashi Takebe
\\
Department of Mathematics\\
Ochanomizu University\\
Otsuka 2-1-1, Bunkyo-ku\\
Tokyo, 112-8610, Japan
}
\date{February 2003; Revised: March 2003}

\maketitle


\begin{abstract}
 Trigonometric degeneration of the Baxter-Belavin elliptic $r$ matrix is
 described by the degeneration of the twisted WZW model on elliptic
 curves. The spaces of conformal blocks and conformal coinvariants of
 the degenerate model are factorised into those of the orbifold WZW
 model.

\end{abstract}

\section{Introduction}
\label{sec:introduction}

The Baxter-Belavin elliptic classical $r$ matrix is a quasi-classical
limit of the famous Baxter's elliptic $R$ matrix.\cite{bax:73} It has a
unique position in the classification of the classical $r$ matrices by
Belavin and Drinfeld.\cite{bel-dri:82}

Etingof introduced the elliptic Knizhnik-Zamolodchikov (KZ) equation as
a system of differential equations satisfied by the twisted trace of
product of vertex operators (intertwining operators of representations)
of affine Lie algebras $\widehat{sl}_N$.\cite{eti:94} This system is
formally obtained by replacing the rational $r$ matrix in the ordinary
KZ equation by the elliptic one. Kuroki and the author constructed the
{\em twisted Wess-Zumino-Witten (WZW) model} on elliptic curves and gave
to Etingof's equations a geometric interpretation as flat
connections.\cite{kur-tak:97}

As is well-known, the elliptic $r$ matrix degenerates to the
trigonometric $r$ matrix when the elliptic modulus $\tau$ tends to
$i\infty$. Hence it is natural to expect that the twisted WZW model
becomes a degenerate model which is related to the trigonometric $r$
matrix when the elliptic curve degenerates to a singular curve.

The purpose of this and the subsequent papers is to show that this is
indeed the case and to analyse the degenerate WZW model. In the present
paper we examine the degenerate WZW model and prove the factorisation
theorem of Tsuchiya-Ueno-Yamada type.\cite{tuy:89} The spaces of
conformal blocks and conformal coinvariants of the twisted WZW model over
the singular curve factorise into those of the orbifold WZW model, which
has been intensively studied in other context. See, e.g.,
\Refcite{orb-cft}. Actually, though based on a different standpoint,
Etingof's work\cite{eti:94} could be considered as a precursor to them.

In the forthcoming papers we study the behaviour of the sheaves of
conformal blocks and conformal coinvariants of the twisted WZW model
around the boundary of the moduli of elliptic curves, especially the
factorisation property of these sheaves. The integrable systems arising
from this WZW model, i.e., the trigonometric Gaudin model and the
trigonometric KZ equations, shall be also analysed.

In this work, we use the formulation of WZW models in terms of the
affine Lie algebras. The reader can translate it into the formulation by
the vertex operator algebras\cite{fre-ben:01,fre-szc:01,nag-tsu:02}
without difficulties.

This paper is organised as follows:
After reviewing the twisted WZW model on elliptic curves in
\secref{sec:twwzw-ell}, we define the trigonometric twisted WZW model
and corresponding WZW model on the orbifold in
\secref{sec:trig,orb-WZW}. The space of conformal coinvariants of the
trigonometric twisted WZW model is decomposed into the ``direct
integral'' of the spaces of conformal coinvariants of the orbifold WZW
model. This general factorisation theorem is proved in
\secref{sec:factorisation1}. If certain conditions are imposed on the
modules inserted to the degenerate curve, a refined factorisation
theorem can be proved. This is done in \secref{sec:factorisation2}. The
simplest non-trivial example of the factorisation is given in
\secref{sec:weyl-mod}. The last section \secref{summary+remarks}
summarises the results and lists several problems.

\subsection*{Notations}

We use the following notations besides other ordinary conventions in
mathematics. 

\begin{itemize}
 \item $N$, $L$: fixed integers. $N \geqq 2$ will be the rank of the Lie
       algebra and $L \geqq 1$ will be the number of the marked points
       on a curve.

 \item Let $X$ be an algebraic variety, $\F$ be a sheaf on $X$ and $P$
       be a point on $X$.

       $\calO_X$, $\Omega_X^1$, $\F_P$, $\m_P$: the structure sheaf of $X$,
       the sheaf of holomorphic 1-forms on $X$, the stalk of $\F$ at
       $P$, the maximal ideal of the local ring $\calO_{X,P}$,
       respectively. 

       When $\F$ is an $\calO_X$-module;

       $\F|_P := \F_P/\m_P\F_P$, $\F^\wedge_P$: fibre of $\F$ at $P$,
       $\m_P$-adic completion of $\F_P$, respectively.

 \item We shall use the same symbol for a vector bundle and for a
       locally free coherent $\calO_X$-module consisting of its local
       holomorphic  sections.

\end{itemize}

\section{Review of twisted WZW model on elliptic curves}
\label{sec:twwzw-ell}

In this section we briefly review the twisted WZW model on elliptic
curves. In the next section the degeneration of this model, the
trigonometric twisted WZW model, shall be given. We follow
\Refcite{kur-tak:97} but make slight changes of normalisation and
notations because we need to use the results of Etingof\cite{eti:94}
later.

We fix an invariant inner product of $\g = sl_N(\C)$ by
\begin{equation}
    (A | B) := \tr (AB) \quad \text{for $A, B \in \g$}.
\label{def:inner-product}
\end{equation}
Define matrices $\beta$ and $\gamma$ by
\begin{equation}
    \beta :=  \begin{pmatrix}
         0 & & 1 & \ \ &       & \ \ & 0 \\
           & & 0 & \ \ &\ddots & \ \ &   \\
           & &   & \ \ &\ddots & \ \ & 1 \\
         1 & &   & \ \ &       & \ \ & 0
         \end{pmatrix},
    \quad
    \gamma := \begin{pmatrix}
         1 &            &        & 0 \\
           & \eps^{-1}  &        &   \\
           &            & \ddots &   \\
         0 &            &        & \eps^{1-N}
         \end{pmatrix},
\label{def:beta,gamma}
\end{equation}
where $\eps = \exp (2\pi i/ N)$.  Then we have $\beta^N=\gamma^N=1$
and $\gamma\beta=\eps\beta\gamma$.

Let $E = E_\tau$ be the elliptic curve with modulus $\tau$: $E_\tau :=
\C/\Z + \tau \Z$. We define the group bundle $\Gtw$ with fibre $G =
SL_N(\C)$ and its associated Lie algebra bundle $\gtw$ with fibre $\g =
sl_N(\C)$ over $E$ by
\begin{align}
    \Gtw :=& (\C \times G)/{\sim},
\label{def:Gtw}
\\
    \gtw :=& (\C \times \g)/{\approx},
\label{def:gtw}
\end{align}
where the equivalence relations $\sim$ and $\approx$ are defined by
\begin{align}
   (z, g) &\sim (z+1, \gamma g \gamma^{-1}) 
          \sim (z+\tau, \beta g \beta^{-1}),
\label{def:eq-rel-Gtw}
\\
   (z, A) &\approx (z+1, \Ad\gamma (A)) 
          \approx (z+\tau, \Ad\beta (A)).
\label{def:eq-rel-gtw}
\end{align}
(Because $\gamma\beta = \eps\beta\gamma$, the group bundle $\Gtw$ is
{\em not} a principal $G$-bundle.\footnote{The bundle $\Gtw$ is a
$PGL_N$ principal bundle. Hurtubise and Markman (for example, see
\Refcite{hur-mar:02}) take this viewpoint, but we prefer to regard it as
a group scheme over the elliptic curve by the reason we stated in \S0 of
\Refcite{kur-tak:97}.}) Let $J_{ab} = \beta^a \gamma^{-b}$, which
satisfies
\begin{equation}
    \Ad \gamma (J_{ab}) = \eps^a J_{ab},\qquad 
    \Ad \beta (J_{ab}) = \eps^b J_{ab}.
\label{Ad(*)Jab}
\end{equation}
Global meromorphic sections of $\gtw$ are linear combinations of $J_{ab}
f(z)$ ($a,b = 0,\dots,N-1$, $(a,b)\neq (0,0)$), where $f(z)$ is a
meromorphic function with quasi-periodicity,
\begin{equation}
    f(z+1) = \eps^a f(z), \qquad f(z+\tau) = \eps^b f(z).
\label{period:gtw:ell}
\end{equation}

The twisted Lie algebra bundle $\gtw$ has a natural connection,
$\nabla_{d/dz} = d/dz$. We define the invariant $\calO_E$-inner product on
 $\gtw$ (regarded as an $\calO_E$-Lie algebra sheaf) by
 \begin{equation}
   (A|B) := \frac{1}{2N} \tr_{\gtw}(\ad A \ad B) \in \calO_E
 \label{def:inner-prod-gtw}
 \end{equation}
 for sections $A$, $B$ of $\gtw$, where the symbol $\ad$ denotes the
 adjoint representation of the $\calO_E$-Lie algebra $\gtw$.  This inner
 product is invariant under the translations with respect to
 the connection $\nabla:\gtw\to\gtw\tensor_{\calO_E}\Omega_E^1$:
 \begin{equation}
   d(A|B) = (\nabla A|B) + (A|\nabla B) \in \Omega_E^1
   \quad
   \text{for $A,B\in\gtw$}.
 \label{compati-conn-innerprod}
 \end{equation}
 Under the trivialisation of $\gtw$ defined by the construction
 \eqref{def:gtw}, the connection $\nabla$ and the inner product
 $(\cdot|\cdot)$ on $\gtw$ coincide with the exterior derivative by $z$
 and the inner product defined by \eqref{def:inner-product}
 respectively.

 For each point $P$ on $E$, we define Lie algebras,
 \begin{equation}
     \g^P := (\gtw \tensor_{\calO_E} \K_E)^{\wedge}_P,
 \label{def:loop}
 \end{equation}
 where $\K_E$ is the sheaf of meromorphic functions on $E$ and
 $(\cdot)_P^{\wedge}$ means the completion of the stalk at $P$ with
 respect to the maximal ideal $\m_P$ of $\calO_P$. In other words, it is
 isomorphic to the loop Lie algebra $\g((z-z_0))$, where $z_0$ is the
 coordinate of $P$, but the isomorphism depends on choices of the
 coordinate $z$ and the trivialisation of $\gtw$. The subspace
 \begin{equation}
     \g^P_+:=(\gtw)_P^{\wedge}\simeqq\g[[z-z_0]]
 \label{def:loop+}
 \end{equation}
 of $\g^P$ is a Lie subalgebra.

Let us fix mutually distinct points $Q_1, \ldots, Q_L$ on $E$ whose
coordinates are $z = z_1, \ldots, z_L$ and put $D:=\{Q_1,\dots,Q_L\}$.
We shall also regard $D$ as a divisor on $E$ (i.e.,
$D=Q_1+\cdots+Q_L$). The Lie algebra $\g^D:=\bigoplus_{i=1}^L\g^{Q_i}$
has a 2-cocycle defined by
\begin{equation}
    \caff(A,B) := \sum_{i=1}^L \caff_{,i}(A_i,B_i),\qquad
    \caff_{,i}(A_i,B_i) := \Res_{Q_i} (\nabla A_i | B_i),
\label{def:cocycle}
\end{equation}
where $A=(A_i)_{i=1}^L, B=(B_i)_{i=1}^L \in \g^D$ and $\Res_{Q_i}$ is
the residue at $Q_i$. (The symbol ``$\caff$'' stands for ``Cocycle
defining the Affine Lie algebra''.)  We denote the central extension of
$\g^D$ with respect to this cocycle by $\hat\g^D$:
\begin{equation*}
    \hat\g^D := \g^D \oplus \C \khat,
\end{equation*}
where $\khat$ is a central element. Explicitly the bracket of $\hat\g^D$ is
represented as
\begin{equation}
  [A, B] = ([A_i, B_i]^\circ)_{i=1}^L \oplus \caff(A,B) \khat
  \quad
  \text{for $A,B\in\g^D$,}
\label{def:aff-alg-str}  
\end{equation}
where $[A_i, B_i]^\circ$ are the natural bracket in $\g^{Q_i}$.  The Lie
algebra $\hat\g^P$ for a point $P$ is nothing but the affine Lie algebra
$\hat\g$ of type $A^{(1)}_{N-1}$ (a central extension of the loop
algebra $\g((t-z))=sl_N\bigl(\C((t-z))\bigr)$).  

The affine Lie algebra $\hat\g^{Q_i}$ can be regarded as a subalgebra of
$\hat\g^D$.  The subalgebra $\g^{Q_i}_+$ of $\g^{Q_i}$ (cf.\
\eqref{def:loop+}) can be also regarded as a subalgebra of
$\hat\g^{Q_i}$ and $\hat\g^D$.

Let $\gout$ be the space of global meromorphic sections of $\gtw$ which are
holomorphic on $E$ except at $D$:
\begin{equation}
    \gout := \Gamma(E, \gtw(\ast D)).
\label{def:gout:ell}
\end{equation}
There is a natural linear map from $\gout$ into $\g^D$ which maps a
meromorphic section of $\gtw$ to its germ at $Q_i$'s. The residue
theorem implies that this linear map is extended to a Lie algebra
injection from $\gout$ into $\hat\g^D$, which allows us to regard
$\gout$ as a subalgebra of $\hat\g^D$.

\begin{definition}
\label{def:CC,CB:ell}
The space of {\em conformal coinvariants} $\CC_E(M)$ and that of {\em
conformal blocks} $\CB_E(M)$ over $E$ associated to
$\hat\g^{Q_i}$-modules $M_i$ with the same level $\khat = k$ are defined
to be the space of coinvariants of $M := \bigotimes_{i=1}^L M_i$ with
respect to $\gout$ and its dual:
\begin{equation}
    \CC_E(M) := M/\gout M, \quad
    \CB_E(M) := \Hom_\C(M/\gout M, \C).
\label{def:cc,cb:ell}
\end{equation}
(In \Refcite{tuy:89}, $\CC_E(M)$ and $\CB_E(M)$ are called the space of
{\em covacua} and that of {\em vacua} respectively.)
\end{definition}

\section{Trigonometric twisted WZW model and orbifold WZW model}
\label{sec:trig,orb-WZW}

In this section we define the degeneration of the twisted WZW model in
the previous section, \secref{sec:twwzw-ell}, i.e., the trigonometric
twisted WZW model on a singular ``curve''. (Exactly speaking, the model
is defined on a ``stack'', whatever it is.) As is shown in the sections
coming later, \secref{sec:factorisation1} and
\secref{sec:factorisation2}, if we desingularise the ``curve'', the
spaces of conformal coinvariants and conformal blocks are factorised
into those of the orbifold WZW model, the definition of which is also
given in this section.

When the modulus $\tau$ of the elliptic curve $E_\tau$ tends to
$i\infty$, one of the cycle of the torus is pinched and we obtain a
singular curve of the ``croissant bread'' shape. Such an object is
constructed simply by identifying the north pole ($\infty$) and the
south pole ($0$) of the Riemann sphere $\Proj^1(\C)$.

Unfortunately, the twisted bundle $\gtw$ over $E_\tau$ does not simply
degenerates to a bundle over such a singular curve. We need to take
several steps:
\begin{enumerate}
 \item Lift the bundle $\gtw$ to a covering space of $E_\tau$ where the
       pull-back of $\gtw$ become trivial. The covering is $N^2$-fold
       and the covering space is an elliptic curve isomorphic to
       $E_\tau$.
 \item Take such degeneration of the covering elliptic curve that it
       becomes the $N^2$-fold covering of the singular curve mentioned
       above. This covering space consists of $N$ copies of $\Proj^1$
       connected to form a ring.
 \item The ``degenerate bundle'' $\gtw$ over the singular curve is
       understood as the trivial bundle over the singular covering space
       constructed above with the cyclic group action together.
\end{enumerate}

\begin{figure}
 \begin{center}
  \includegraphics{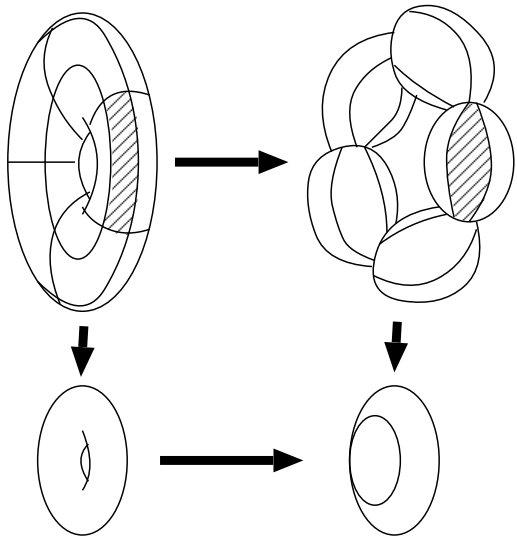}
 \end{center}
\vspace*{8pt}
\caption{Degeneration of an elliptic curve and its $N^2$-covering.}
\end{figure}

A ``manifold with local finite group action'' is called an {\em
orbifold}, but in our case the underlying space has a singularity, hence
it should be regarded as a {\em stack} in algebraic geometry. It is not
our task (and certainly beyond the reach of the author's ability) to
construct a general theory of the WZW model over stacks. We shall define
the twisted WZW model over the degeneration of the elliptic curve below
by down-to-earth terminology only.

Let us denote the standard coordinate of $\Proj^1(\C)$ by $t$. The
cyclic group $C_N = \Z/N\Z$ acts as $t \mapsto \eps^a t$ ($a \in C_N$)
and the quotient $E_\orb = \Proj^1/C_N$ is an ordinary orbifold. The
degenerate elliptic curve $E_\trig$ is defined by identifying $0$ and
$\infty$ (i.e., the fixed points of the $C_N$ action) of $E_\orb$.

Not defining what the twisted bundle $\g^\trig$ over $E_\trig$ and the
twisted bundle $\g^\orb$ over $E_\orb$ are, we directly define what
their meromorphic sections are in the following way. (In fact,
specifying sections is equivalent to defining bundles in algebraic
geometry.)

\begin{definition}
\label{def:gtrig,gorb}
A $\g$-valued meromorphic function $f(t)$ on $\Proj^1$ is a meromorphic
section of $\g^\orb$ if and only if it satisfies
\begin{equation}
   f(\eps t) = \Ad\gamma (f(t)).
\label{period:gorb}
\end{equation}
The same $f(t)$ is a section of $\g^\trig$ if and only if it satisfies
\begin{equation}
   f(\infty) = \Ad\beta (f(0)),
\label{period:gtrig}
\end{equation}
in addition to \eqref{period:gorb}.
\end{definition}

The relations \eqref{period:gtrig} and \eqref{period:gorb} are
degenerate form of \eqref{def:eq-rel-gtw}. 

Hereafter we fix mutually distinct points $Q_1, \ldots, Q_L$ on the
orbifold $E_\orb$ with coordinates $t = t_1, \ldots, t_L$ (modulo
$C_N$-action). The algebras $\g^{\trig,Q_i}$, $\g^{\trig,Q_i}_+$,
$\g^{\orb,Q_i}$ and $\g^{\orb,Q_i}_+$ are defined from $\g^\trig$ and
$\g^\orb$ in the same way as $\g^{Q_i}$ and $\g^{Q_i}_+$ for the
elliptic curve, \eqref{def:loop} and \eqref{def:loop+}. When $Q_i \neq
0, \infty$ ($i=1,...,L$), they are isomorphic to $\g^{Q_i}$ and
$\g^{Q_i}_+$ since the definition is local.

For the orbifold WZW model, we need similar Lie algebras $\g^{(0)}$,
$\g^{(0)}_+$, $\g^{(\infty)}$, $\g^{(\infty)}_+$ sitting at the singular
points, $0$ and $\infty$. Formally their definition is the same as
\eqref{def:loop} and \eqref{def:loop+} with $P$ replaced by $0$ and
$\infty$ and $E$ by $E_\orb$. The mode expansion of $\g^{(0)}$ has a
different form from the usual loop algebra, namely,
\begin{equation}
    \g^{(0)} = \bigoplus_{\substack{a,b=0,\dots,N-1\\ (a,b)\ne(0,0)}}
    \bigoplus_{m\in\Z} 
    \C J_{a,b} \tensor t^{a+mN}.
\label{g(0)}
\end{equation}
The algebra $\g^{(\infty)}$ has the same form.

The central extension of $\g^{Q_i}$, $\hat\g^{(Q_i)}$, is defined in the
same way as in \secref{sec:twwzw-ell} by using the cocycle
\eqref{def:cocycle}. The central extension of $\g^{(0)}$ and
$\g^{(\infty)}$ are defined essentially in the same manner.
Difference is that we have to change the
normalisation of the central element $\hat k$:
\begin{equation}
    [A, B] = [A, B]^\circ \oplus \caff_{,0}(A,B) \khat,
    \quad
    \caff_{,0}(A,B):= \frac{1}{N}\Res_{t=0} (\nabla A|B)
\label{def:aff-alg-str:0}
\end{equation}
for $A,B\in\g^{(0)}$, where $[A, B]^\circ$ is the natural loop-algebra
bracket in $\g^{{(0)}}$.  The origin of the factor $1/N$ shall be
explained soon later. The definition of $\hat\g^{(\infty)}$ is exactly
the same. We often denote $\g^{(0)}$ and $\g^{(\infty)}$ by
$\g^{(0)}_\gamma$ and $\g^{(\infty)}_\gamma$, following \S1 of
\Refcite{eti:94}, where $\hat\g_\gamma$ is shown to be isomorphic to the
ordinary affine Lie algebra $\hat\g$. The suffix $\gamma$ is put to
emphasise the twisting of the algebra by $\gamma$, \eqref{period:gorb}.

Let $D$ be the divisor $Q_1+ \cdots + Q_L$. The definition of the loop
algebra $\g^D$ and its central extension does not change so far as $D$
does not contain $0$ and $\infty$. When $D$ contains $0$ or $\infty$ or
both, then the change of the normalisation \eqref{def:aff-alg-str:0} of
the cocycle at $0$ and $\infty$ should be taken into account. (When we
talk about $E_\trig$ or $\g^\trig$, it is always assumed that $D$
contains neither $0$ nor $\infty$.)

Instead of the subalgebra $\gout$, \eqref{def:gout:ell}, we define the
following subalgebras of $\hat\g^D$:
\begin{align}
    \gout^\trig &:= 
    \{ f(t): \text{global meromorphic section of\ }\g^\trig
             \text{\ with poles in\ }D\}.
\label{def:gout:trig}
\\
    \gout^\orb &:= 
    \{ f(t): \text{global meromorphic section of\ }\g^\orb 
             \text{\ with poles in\ }D\},
\label{def:gout:orb}
\end{align}
Their Lie algebra structures are defined by the pointwise Lie bracket of
$\g$-valued functions. It is not trivial that they are subalgebras of
$\hat\g^D$, since the residue theorem cannot be directly applied. For
$f(t), g(t) \in \gout^\trig$, it is easy to see that
\begin{equation*}
    \caff(f,g) = \sum_{i=1}^L \Res_{Q_i} (\nabla f(t) | g(t))
    =
    \frac{1}{N} \sum_{a=0}^{N-1} 
    \sum_{i=1}^L \Res_{\eps^a Q_i} (\nabla f(t) | g(t)) = 0,
\end{equation*}
where $\eps^a Q_i$ denotes the action of $a\in C_N$ on the point
$Q_i$. Since $\{ \eps^a Q_i \mid a \in C_N, i=1,...,L\}$ contains all
poles of $(\nabla f(t) | g(t))$ on $\Proj^1$, we can apply the residue
theorem here. This proves that $\gout^\trig$ defined by
\eqref{def:gout:trig} is a Lie subalgebra of $\hat\g^D$. In the case of
$\gout^\orb$, when $D$ does not contain $0$ nor $\infty$, the proof is
the same. When $D$ contains, e.g., $0$, the central extension should be
$\hat k$ times
\begin{equation*}
 \begin{split}
    \caff(f,g) 
    =& 
    \caff_{,0}(f(t),g(t)) + \sum_{i=1}^L \caff_{,i} (f(t) | g(t))
\\
    =&
    \frac{1}{N}\left(\Res_{0} (\nabla f(t) | g(t)) +
    \sum_{a=0}^{N-1} 
    \sum_{i=1}^L \Res_{\eps^a Q_i} (\nabla f(t) | g(t)) \right)
    = 0.
 \end{split}
\end{equation*}
The last equality is the consequence of the residue theorem applied to
$(\nabla f(t)|g(t))$ over $\Proj^1$. Hence, also in this case
$\gout^\orb$ is a subalgebra of $\hat\g^D$. Other cases (when $D$
contains $\infty$) are proved in the same way. This is the reason why we
have to put $1/N$ in \eqref{def:aff-alg-str:0}, which appeared first in
\Refcite{eti:94}, Lemma 1.1, where Etingof wrote down explicit formulae
of the isomorphism of $\hat\g_\gamma$ and the ordinary affine Lie
algebra $\hat\g$.

We are now ready to define the spaces of conformal coinvariants and
conformal blocks for $E_\trig$ and $E_\orb$.

\begin{definition}
\label{def:CC,CB:trig,orb}
(i)
 The space of {\em conformal coinvariants} $\CC_\trig(M)$ and that of
 {\em conformal blocks} $\CB_\trig(M)$ over $E_\trig$ associated to
 $\hat\g^{Q_i}$-modules $M_i$ with the same level $\khat = k$ are
 defined to be the space of coinvariants of $M := \bigotimes_{i=1}^L
 M_i$ with respect to $\gout^\trig$ and its dual:
\begin{equation}
    \CC_\trig(M) := M/\gout^\trig M, \quad
    \CB_\trig(M) := \Hom_\C(M/\gout^\trig M, \C).
\label{def:cc,cb:trig}
\end{equation}

(ii)
 The space of {\em conformal coinvariants} $\CC_\orb(M)$ and that of
 {\em conformal blocks} $\CB_\orb(M)$ over $E_\orb$ associated to
 $\hat\g^{Q_i}$-modules $M_i$ with the same level $\khat = k$ are
 defined to be the space of coinvariants of $M := \bigotimes_{i=1}^L
 M_i$ with respect to $\gout^\orb$ and its dual:
\begin{equation}
    \CC_\orb(M) := M/\gout^\orb M, \quad
    \CB_\orb(M) := \Hom_\C(M/\gout^\orb M, \C).
\label{def:cc,cb:orb}
\end{equation}

\end{definition}

In this paper we consider only {\em smooth modules} of affine Lie
algebras $\hat\g^{Q_i}$, $\hat\g^{(0)}$, etc.: Let $\hat\g = \g\tensor
\C[\xi,\xi^{-1}] \oplus \C\hat k$ be an affine Lie algebra. (We identify
$\hat\g^{Q_i}$, $\hat\g^{(0)}$ and so on with this algebra by means of
suitable isomorphisms.) A $\hat\g$-module $M$ is called {\em smooth} if
for any $v\in M$, there exists $m\geqq 0$ such that, for $A_1, \dots,
A_\nu\in\g$, $m_1,\dots,m_\nu\geqq0$, and $\nu=0,1,2,\ldots$,
\begin{equation}
  (A_1\tensor\xi^{m_1})
  \cdots
  (A_\nu\tensor\xi^{m_\nu})
  v_i=0
  \quad \text{if $m_1+\dots+m_\nu\geqq m$},
\label{smoothness}
\end{equation}
where $\xi$ is the loop parameter of $\hat\g = \g((\xi)) \oplus \C \hat
k$. Not only Laurent polynomials but also Laurent series in
$\g\tensor\C[\xi^{-1}][[\xi]]$ can act on smooth modules. Therefore, if
$M_i$ ($i=1,\dots,L$) are smooth $\hat\g^{Q_i}$ modules, the actions of
$\gout^\trig$ and $\gout^\orb$ on $M_1 \tensor \cdots \tensor M_L$ are
well-defined.

Concrete examples of spaces of conformal coinvariants shall be given in
\secref{sec:weyl-mod}. The relation of $\CC_\trig$ and $\CC_\orb$ is our
main subject in this paper and studied in the next sections,
\S\S\ref{sec:factorisation1}--\ref{sec:factorisation2}. How $\CC_E$ for
elliptic curves degenerates to these spaces at the boundary of the
moduli space of the elliptic curves shall be clear in the forthcoming
papers.

\section{Factorisation (1)}
\label{sec:factorisation1}

In this section the general factorisation theorem of the space of
conformal coinvariants $\CC_\trig$ is shown. The result is not so sharp
as that in the next section, but the idea of the factorisation becomes
clear.

To state the theorem, we need the notion of the Verma module of
$\hat\g_\gamma$, i.e., $\hat\g^{(0)}$ and $\hat\g^{(\infty)}$. Recall
that $\hat\g^{(0)}$ has a natural triangular decomposition:
\begin{equation}
    \hat\g^{(0)} 
    = \n^{(0)}_+ \oplus \h^{(0)} \oplus \n^{(0)}_- \oplus \C \hat k,
\label{decomp:g(0)}
\end{equation}
where $\n^{(0)}_+$ consists of $\g$-valued power series in $t$ with
strictly positive powers while $\n^{(0)}_-$ consists of series in $t$
with negative powers. (See \eqref{g(0)}.) $\h^{(0)}$ is the same as the
Cartan subalgebra of $\g$, which is spanned by $J_{0,b}$
($b=1,...,N-1$). Since $\h^{(0)}$ is canonically identified with the
Cartan subalgebra $\h$ of $\g$, we often drop the superscript $(0)$. The
algebra $\hat\g^{(\infty)}$ has the same decomposition but with opposite
powers of $t$: $\n^{(\infty)}_+$ is the subspace of $\g^{(\infty)}$ with
negative powers of $t$ and $\n^{(\infty)}_-$ is that with negative
powers of $t$, because $\n^{(\infty)}_+$ is the space of $\g$ valued
functions which vanish at $\infty$. To avoid repeating the same
statements for $\hat\g^{(0)}$ and $\hat\g^{(\infty)}$, we identify them
in the canonical way and denote the above decomposition as
\begin{equation}
    \hat\g_\gamma 
    = 
    \n_{\gamma,+} \oplus \h_\gamma \oplus \n_{\gamma,-} \oplus \C \hat k,
\label{decomp:g-gamma}
\end{equation}
and define $\hat\borel_{\gamma} := \n_{\gamma,+} \oplus \h_\gamma \oplus
\C \hat k$. 

\begin{definition}
\label{def:verma}
(i)
 Let $\lambda$ is an element of $\h_\gamma^* = \Hom_\C(\h_\gamma,\C)$
 and $k$ is an arbitrary complex number.  The {\em Verma} module of
 $\hat\g_\gamma$ of highest weight $\lambda$ and level $k$ is the
 $\hat\g_\gamma$-module
\begin{equation}
    M_\lambda := U\hat\g_\gamma / J_{\gamma,\lambda,k},
\label{verma:quot}
\end{equation}
 where $J_{\gamma,\lambda,k}$ is the left ideal generated by
 $\n_{\gamma,+}$, $H - \lambda(H)$ for $H \in \h_\gamma$ and $\hat k -
 k$. The cyclic vector $1 \bmod J_{\gamma,\lambda,k}$ is denoted by
 $1_\lambda$. Equivalently $M_\lambda$ is defined as the
 $\hat\g_\gamma$-module induced from the one-dimensional
 $\hat\borel_\gamma$-module:
\begin{equation}
    M_\lambda := \Ind_{\hat\borel_\gamma}^{\hat\g} \C 1_\lambda,
\label{verma:ind}
\end{equation}
 where $1_\lambda$ satisfies $\n_{\gamma,+} 1_\lambda = 0$, $H 1_\lambda
 = \lambda(H) 1_\lambda$ and $\hat k 1_\lambda = k 1_\lambda$. (See
 \Refcite{eti:94}, (1.3).)

(ii)
 The $\hat\g_\gamma$-module
\begin{equation}
    M(0):= U\hat\g_\gamma/J(0),
\label{univerma}
\end{equation}
 where $J(0) = U\hat\g_\gamma\, \n_{\gamma,+} + U\hat\g_\gamma\, (\hat k
 - k)$, is called the {\em universal Verma} module of level
 $k$.\cite{kas:85} The cyclic vector $1 \bmod J(0)$ is denoted by
 $\univac$. Equivalently $M(0)$ is defined as the $\hat\g_\gamma$-module
 induced from the one-dimensional $(\n_{\gamma,+}\oplus\C\hat
 k)$-module:
\begin{equation}
    M(0) := \Ind_{\n_{\gamma,+}\oplus\C\hat k}^{\hat\g} \C \univac,
\label{univerma:ind}
\end{equation}
 where $\univac$ satisfies $\n_{\gamma,+} \univac = 0$ and $\hat k
 \univac = k \univac$.
\end{definition}

Since $[\h_\gamma, J_{\gamma,\lambda,k}] \subset J_{\gamma,\lambda,k}$
and $[\h_\gamma,J(0)] \subset J(0)$, we can define right
$\h_\gamma$-actions on $M_\lambda$ and $M(0)$ by
\begin{equation}
     [x] \mapsto [x H],
\label{right-h-action:verma}
\end{equation}
for $H \in \h_\gamma$ and $x \in U\hat\g_\gamma$, where $[x]$ denotes
the equivalence class of $x$ modulo $J_{\gamma,\lambda,k}$ or $J(0)$.

The general factorisation theorem is as follows. Hereafter we always
assume that $\{Q_1, \dots, Q_L\}$ do not contain neither $0$ nor
$\infty$. Hence $D = Q_1 + \cdots + Q_L$ is a divisor of $E_\trig$ as
well as that of $E_\orb$.

\begin{theorem}
\label{thm:general-factorisation}
 Let $M_i$ be a smooth $\hat\g^{Q_i}$-module of level $k$ for
 $i=1,\dots,L$ and $M$ be the tensor product $M_1 \tensor \cdots \tensor
 M_L$ endowed with the $\hat\g^D$-module structure. Then,
\begin{equation}
 \begin{split}
    \CC_\trig (M) 
    &\simeqq
    \CC_\orb(M(0)_0\tensor M \tensor M(0)_\infty)/ \CC_\orb\cdot\h,
\\
    \CC_\orb\cdot\h 
    &:=
    \CC_\orb(M(0)_0\tensor M \tensor M(0)_\infty) \rho_{1,\beta}(\h).
 \end{split}
\label{gen-fact}
\end{equation}
 Here $M(0)_0$ and $M(0)_\infty$ are the universal Verma modules defined
 above for $\hat\g_\gamma = \hat\g^{(0)}$ and for $\hat\g_\gamma =
 \hat\g^{(\infty)}$ respectively and $\rho_{1,\beta}$ is the right
 action of $\h = \h^{(0)} = \h^{(\infty)}$ on $M(0)_0\tensor M \tensor
 M(0)_\infty$ defined by
\begin{equation}
    (v_0 \tensor v \tensor v_\infty) \rho_{1,\beta}(H)
    :=
    v_0 H \tensor v \tensor v_\infty +
    v_0   \tensor v \tensor v_\infty \Ad(\beta) H
\label{right-h-action:cc}
\end{equation}
 for $H \in \h$, $v_0 \in M(0)_0$, $v \in M$ and $v_\infty \in
 M(0)_\infty$. Note that the left action of $\gout^\orb$ and the right
 action of $\h$ on $M(0)_0 \tensor M \tensor M(0)_\infty$ commute with
 each other and thus $\rho_{1,\beta}$ descends to the action on
 $\CC_\orb$. 
\end{theorem}

This statement is rigorous but not very illuminating, so we rephrase it
in a little bit incorrect but more impressive way. By virtue of the
Poincar\'e-Birkhoff-Witt theorem, $M_\lambda$ and $M(0)$ are isomorphic
to $U\n_{\gamma,-}$ and $U\n_{\gamma,-} \tensor U\h_\gamma$ respectively
as linear spaces. Being commutative, the universal enveloping algebra
$U\h_\gamma$ is a polynomial ring over the dual space $\h_\gamma^* =
\Hom_\C(\h_\gamma,\C)$. Hence $M(0)$ is spanned by $x\tensor
f(\lambda)$, where $x \in U\g_{\gamma,-}$ and $f(\lambda)$ is a
polynomial in $\lambda \in \h_{\gamma}^*$. Evaluating $x \tensor
f(\lambda)$ at $\lambda \in \h_\gamma^*$, we obtain an element of
$M_\lambda$. In this sense, $M(0)$ may be approximately regarded as a
direct integral:
\begin{equation}
    M(0)\  \text{``$=$''} 
    \int_{\lambda\in\h_\gamma^*} M_\lambda\, d\lambda.
\label{M(0)=int(Mlambda)}
\end{equation}
(This is not correct since the direct integral should be spanned by
$L^2$-functions instead of polynomials.) Using this (approximate)
identification, we can rewrite $\CC_\orb(M(0)_0 \tensor M \tensor
M(0)_\infty)$ as
\begin{equation}
    \CC_\orb(M(0)_0 \tensor M \tensor M(0)_\infty) 
    = \iint_{\h^{*} \times \h^{*}}
    \CC_\orb(M_\lambda \tensor M \tensor M_{\mu})\, 
    d\lambda \, d\mu.
\end{equation}
The statement of \thmref{thm:general-factorisation} becomes
\begin{equation}
     \CC_\trig(M)
     \simeqq
     \int_{\h}
     \CC_\orb(M_{\lambda}\tensor M \tensor M_{\beta^*(\lambda)})
     \, d\lambda,
\label{factorisation:integral}
\end{equation}
where $\beta^*(\lambda) = - \lambda \circ \Ad(\beta^{-1})$. The proof of
\eqref{factorisation:integral} (from \thmref{thm:general-factorisation})
is the same as that of \eqref{fact:finite:component}.

The isomorphism \eqref{factorisation:integral} is an analogue of the
factorisation theorem of ordinary WZW
models.\cite{tuy:89,uen:95,shi-uen:99,nag-tsu:02} In fact, under certain
finiteness conditions for $M_i$'s, almost all
$\CC_\orb(M_{\lambda}\tensor M \tensor M_{\beta^*(\lambda)})$ vanish and
the direct integral is replaced by finite direct sum as we will prove in
\secref{sec:factorisation2}.

\begin{proof}
 The idea of the proof of \thmref{thm:general-factorisation} is
 essentially the same as that of \Refcite{tuy:89}. See also
 \Refcite{uen:95,shi-uen:99}. We make use of techniques in
 \Refcite{nag-tsu:02}.

 To begin with, observe that there is an exact sequence,
\begin{equation}
    0 \to 
    \gout^\trig M/ \gout^0 M \to
    M/ \gout^0 M \to 
    \CC_\trig(M) \to 0,
\label{ex-seq:CCtrig}
\end{equation}
 where 
\begin{multline}
    \gout^0 = \{ f(t): \text{meromorphic section of\ }\g^\orb
         \text{\ with poles in\ }D, f(0) = f(\infty) = 0\} 
\\
     \subset \gout^\trig.
\label{def:gout0}
\end{multline}
 Using this exact sequence, we prove the theorem in two steps. In the
 first step, we show that a linear map $\varphi(v \bmod \gout^0 M) =
 [\univac_0 \tensor v \tensor \univac_\infty] \in \CC_\orb(M(0)_0
 \tensor M \tensor M(0)_\infty)$ induces an isomorphism
\begin{equation}
    \varphi: M/\gout^0 M \overset{\sim}{\to}
    \CC_\orb(M(0)_0 \tensor M \tensor M(0)_\infty).
\label{phi:isom}
\end{equation}
 In the second step we show that 
\begin{equation}
    \varphi(\gout^\trig M/\gout^0 M)
    \simeqq
    \CC_\orb\cdot\h,
\label{kernel:M/g0->CCtrig:0}
\end{equation}
 (cf.\ \eqref{gen-fact} for the definition of $\CC_\orb\cdot\h$) or,
 equivalently, 
\begin{equation}
    \univac_0 \tensor \gout^\trig M \tensor \univac_\infty
    \equiv 
    (M(0)_0 \tensor M \tensor M(0)_\infty) \rho_{1,\beta}(\h)
    \mod \gout^\orb.
\label{kernel:M/g0->CCtrig:1}
\end{equation}
 Hereafter $\mod \gout^\orb$ denotes $\mod \gout^\orb(M(0)_0 \tensor M
 \tensor M(0)_\infty)$.
 Combining these two statements with the exact sequence
 \eqref{ex-seq:CCtrig}, we obtain the proof of the theorem.

 To prove that $\varphi$ in \eqref{phi:isom} is an isomorphism, we use a
 lemma by Kazhdan and Lusztig.\cite{kaz-lus:93}

\begin{lemma}
\label{lem:KL}
 Let $\germa$ be a Lie algebra over $\C$ and $\germa_1$, $\germa_2$ be
 its Lie subalgebra such that $\germa = \germa_1 + \germa_2$. Assume
 that an $\germa_2$-module $V$ is given and put $W :=
 \Ind_{\germa_2}^\germa V$. Then the canonical homomorphism $V \owns v
 \mapsto v \in W$ induces an isomorphism
\begin{equation}
     V/(\germa_1 \cap \germa_2) V
     \overset{\sim}{\to}
     W/\germa_1 W.
\label{isom:KL}
\end{equation}
\end{lemma}
 See \S A.7 of \Refcite{kaz-lus:93} for the proof.

 We apply this lemma to
\begin{equation}
    \germa_1 = \gout^\orb, \qquad
    \germa_2 = \n^{(0)}_+ \oplus \g^D \oplus \n^{(\infty)}_+
               \oplus \C \hat k, \qquad
    V = \C \univac_0 \tensor M \tensor \C \univac_\infty,
\end{equation}
 where $\univac_0$ and $\univac_\infty$ are the canonical cyclic vectors
 of $M(0)_0$ and $M(0)_\infty$. It is easy to see from
 \eqref{decomp:g(0)} or \eqref{decomp:g-gamma} that the Lie algebra
 $\germa = \germa_1 + \germa_2$ is equal to $\g^{(0)} \oplus \g^D \oplus
 \g^{(\infty)} \oplus \C \hat k$. Since $\germa_1 \cap \germa_2 =
 \gout^0$ acts trivially on the first and the last factor of
 $M(0)_0 \tensor M \tensor M(0)_\infty$, we have
\begin{equation}
    V/(\germa_1 \cap \germa_2) V
    \simeqq
    \C \univac_0 \tensor M \tensor \C \univac_\infty /
    \gout^0 \C (\univac_0 \tensor M \tensor \C \univac_\infty)
    = M/\gout^0 M.
\label{Va=M/gout0M}
\end{equation}
 On the other hand, 
\begin{equation}
    W = 
    \Ind_{\n^{(0)}_+ \oplus \g^D \oplus \n^{(\infty)}_+ \oplus \C \hat k}
        ^{\g^{(0)} \oplus \g^D \oplus \g^{(\infty)} \oplus \C \hat k}
    \C \univac_0 \tensor M \tensor \C \univac_\infty
    = M(0)_0 \tensor M \tensor M(0)_\infty.
\end{equation}
 Hence $W/\germa_1 W = M(0)_0 \tensor M \tensor
 M(0)_\infty/\gout^\orb(M(0)_0 \tensor M \tensor M(0)_\infty)$, which
 proves \eqref{phi:isom} by \lemref{lem:KL}.

 Next we prove \eqref{kernel:M/g0->CCtrig:1}. Namely, we prove the
 following two claims:
\begin{gather}
    \univac_0 \tensor \gout^\trig M \tensor \univac_\infty \subset
    (M(0)_0 \tensor M \tensor M(0)_\infty) \rho_{1,\beta}(\h) +
    \gout^\orb (M(0)_0 \tensor M \tensor M(0)_\infty).
\label{kernel:M/g0->CCtrig:2}
\\
    \univac_0 \tensor \gout^\trig M \tensor \univac_\infty +
    \gout^\orb (M(0)_0 \tensor M \tensor M(0)_\infty) \supset 
    (M(0)_0 \tensor M \tensor M(0)_\infty) \rho_{1,\beta}(\h).
\label{kernel:M/g0->CCtrig:3}
\end{gather}

 To prove \eqref{kernel:M/g0->CCtrig:2}, assume $f(t) \in \gout^\trig$
 and $v \in M$. We show that $\univac_0 \tensor f(t) v \tensor
 \univac_\infty$ belongs to $(M(0)_0 \tensor M \tensor M(0)_\infty)
 \rho_{1,\beta}(\h) \mod \gout^\orb$. Since $f(t) \in \gout^\trig \subset
 \gout^\orb$, we have
\begin{equation}
    \univac_0 \tensor f(t) v \tensor \univac_\infty 
    \equiv 
    - f(0) \univac_0 \tensor v \tensor \univac_\infty
    - \univac_0 \tensor v \tensor f(\infty) \univac_\infty
    \mod{\gout^\orb}.
\label{1*fv*1}
\end{equation}
 The conditions \eqref{period:gorb} and \eqref{period:gtrig} imply that
 $f(\infty) = \Ad\beta (f(0)) \in \h$. Hence the right hand side of
 \eqref{1*fv*1} belongs to $(M(0)_0 \tensor M \tensor M(0)_\infty)
 \rho_{1,\beta}(\h)$ which proves \eqref{kernel:M/g0->CCtrig:2}.

 Conversely, assume $v_0 \in M(0)_0$, $v \in M$, $v_\infty \in
 M(0)_\infty$ and $H \in \h$. We show that $(v_0 \tensor v \tensor
 v_\infty) \rho_{1,\beta}(H)$ belongs to $\univac_0 \tensor \gout^\trig
 M \tensor \univac_\infty \mod{\gout^\orb}$ to prove
 \eqref{kernel:M/g0->CCtrig:3}. Thanks to the fact we proved in the
 first step of this proof of \thmref{thm:general-factorisation}, there
 exists $\tilde v \in M$ such that
\begin{equation}
    v_0 \tensor v \tensor v_\infty \equiv 
    \univac_0 \tensor \tilde v  \tensor \univac_\infty
    \mod{\gout^\orb}.
\label{def:tilde-v}
\end{equation}
 Hence we have
\begin{equation}
 \begin{split}
    (v_0 \tensor v \tensor v_\infty) \rho_{1,\beta}(H)
    &\equiv
    (\univac_0 \tensor \tilde v \tensor \univac_\infty) \rho_{1,\beta}(H) 
\\
    &=
    H \tensor \tilde v \tensor \univac_\infty
    +
    \univac_0 \tensor \tilde v \tensor \Ad\beta(H) 
\\
    &\equiv
    - \univac_0 \tensor h(t) \tilde v  \tensor \univac_\infty
 \end{split}
\label{1*h(t)v*1}
\end{equation}
 modulo $\gout^\orb (M(0)_0 \tensor M \tensor M(0)_\infty)$. Here $h(t)$
 is an element of $\gout^\trig$ satisfying $h(0) = H + (\text{regular
 function at\ }t=0)$ and $h(\infty) = \Ad\beta(H) + (\text{regular
 function at\ }t=\infty)$. When $H$ is $J_{0,b}$ (cf.\
 \eqref{Ad(*)Jab}), such $h(t)$ can be written down explicitly:
\begin{equation}
    h_b(t) := 
    \frac{\eps^b t^N - t_i^N}{t^N - t_i^N} J_{0,b},
\label{def:h(t)}
\end{equation}
 where $t_i = \exp(2\pi i z_i/N)$ is the coordinate of $Q_i$ for a
 certain $i$ ($1\leqq i \leqq N$). For general $H$, $h(t)$ can be
 constructed as a linear combination of $h_b(t)$. The last expression in
 \eqref{1*h(t)v*1} is an element of $\univac_0 \tensor \gout^\trig M
 \tensor \univac_\infty$, which proves
 \eqref{kernel:M/g0->CCtrig:3}. Thus \eqref{kernel:M/g0->CCtrig:1} as
 well as \thmref{thm:general-factorisation} is proved.
\end{proof}

\section{Factorisation (2)}
\label{sec:factorisation2}

The statement of the previous section can be refined for special
modules. In fact, the direct integral in \eqref{factorisation:integral}
is replaced by finite direct sum when we impose following finiteness
conditions on $M_i$'s: for each $i=1,\dots,L$ there exists a finite
dimensional $\g$-module $V_i$ which generates $M_i$ over $\hat\g^{Q_i}$
and $\hat\g^{Q_i}_+ V_i \subset V_i$. Throughout this section we assume
this condition, which automatically implies that $M_i$'s are smooth.

Let us denote $V=V_1 \tensor \cdots \tensor V_L$ and the set of
$\g$-weights of $V$ by $\wt(V)$. For $\lambda \in \h^*$, we define 
\begin{equation}
 \begin{split}
    \tilde \lambda  &:= - \lambda \circ (1- \Ad\beta^{-1})^{-1}
    = \lambda \circ (1- \Ad\beta)^{-1} \circ \Ad\beta,
\\
    \tilde \lambda' &:= - \lambda \circ (1- \Ad\beta)^{-1}.
 \end{split}
\label{def:tilde-lambda}
\end{equation}
They are well-defined since $1 - \Ad\beta^{\pm1}$ is invertible on $\h =
\bigoplus_{b=1}^{N-1} \C J_{0,b}$.

\begin{theorem}
\label{thm:factorisation:finite}
 Under the above assumption on $M_i$, $\CC_\trig(M)$ is finite
 dimensional and there is a linear isomorphism
\begin{equation}
    \CC_\trig(M)
    \simeqq
    \bigoplus_{\lambda \in \wt(V)}
    \CC_\orb(M_{\tilde\lambda} \tensor M \tensor M_{\tilde\lambda'}),
\label{factorisation:finite}
\end{equation}
 where $M_{\tilde\lambda}$ and $M_{\tilde\lambda'}$ are inserted at $0$
 and $\infty$ of $E_\orb$ respectively.
\end{theorem}

\begin{proof}
 For simplicity of notation, we denote $J_{0,b}$ by $H_b$
 ($b=1,\dots,N-1$). As mentioned before, they form a basis of $\h$. We
 also introduce ideals $I_\h$ and $I'_\h$ of the polynomial ring $U\h$
 as follows:
\begin{equation}
 \begin{gathered}
    I_\h  = \bigcap_{\lambda\in\wt(V)} I_\lambda,\qquad
    I'_\h = \bigcap_{\lambda\in\wt(V)} I'_\lambda,
\\
    I_\lambda := \text{ideal of\ }U\h \text{\ generated by\ }
    H_b - \tilde\lambda(H_b) (b=1,\dots,N-1).
\\
    I'_\lambda := \text{ideal of\ }U\h \text{\ generated by\ }
    H_b - \tilde\lambda'(H_b) (b=1,\dots,N-1).
 \end{gathered}
\label{def:Ih,I'h}
\end{equation}

 The main part of the proof of the theorem is to show the following
 inclusion:
\begin{equation}
    I \tensor M \tensor M(0)_\infty + M(0)_0 \tensor M \tensor I'
    \subset K,
\label{kernel-inclusion}
\end{equation}
 where 
\begin{equation}
    I  = M(0)_0 I_\h, \qquad I' = M(0)_\infty I'_\h, 
\label{def:I,I'}
\end{equation}
 and $K$ is the kernel of the canonical surjection $M(0)_0 \tensor
 M \tensor M(0)_\infty \to \CC_\orb(M(0)_0\tensor M \tensor
 M(0)_\infty)/ \CC_\orb\cdot\h$, namely,
\begin{equation}
    K :=
    \gout^\orb (M(0)_0 \tensor M \tensor M(0)_\infty)
    + (M(0)_0 \tensor M \tensor M(0)_\infty) \rho_{1,\beta}(\h).
\label{def:K}
\end{equation}
 When the inclusion \eqref{kernel-inclusion} is proved, we shall see
 that we may replace $\CC_\orb(M(0)_0 \tensor M \tensor M(0)_\infty)$ in
 \eqref{gen-fact} by $\CC_\orb(\bar M(0)_0 \tensor M \tensor \bar
 M(0)_\infty)$, where $\bar M(0)_0 = M(0)_0 / I$ and $\bar M(0)_\infty =
 M(0)_\infty/I'$ which are finite direct sums of Verma modules. The
 space of conformal coinvariants $\CC_\orb(\bar M(0)_0 \tensor M \tensor
 \bar M(0)_\infty)$ decomposes into finite direct sum of
 $\CC_\orb(M_{\tilde\lambda} \tensor M \tensor M_{\tilde\lambda'})$ and,
 factoring by $\CC_\orb\cdot \h$, we shall obtain
 \eqref{factorisation:finite}. The proof of the first statement, $\dim
 \CC_\trig(M)< \infty$, is given after \lemref{lem:decomp-gD}.

 Let us show that
\begin{equation}
    I \tensor M \tensor M(0)_\infty \subset K,
\label{I*M*M<K}
\end{equation}
 namely, for $v_0 \tensor v \tensor v_\infty \in M(0)_0 \tensor M
 \tensor M(0)_\infty$ and $p(H) \in I_\h$, $v_0 p(H) \tensor v \tensor
 v_\infty$ belongs to $K$. Thanks to the existence of $\tilde v$
 satisfying \eqref{def:tilde-v}, we may assume $v_0 = \univac_0$ and
 $v_\infty = \univac_\infty$, because $\gout^\orb(M(0)_0 \tensor M
 \tensor M(0)_\infty) (p(H) \tensor 1 \tensor 1) \subset K$. Hence the
 problem reduces to show
\begin{equation}
    p(H) \tensor v \tensor \univac_\infty \in K.
\label{p(H)*v*1<K}
\end{equation}
 By the following lemma, we have only to show \eqref{p(H)*v*1<K} for $v
 \in V \subset M$.

\begin{lemma}
\label{lem:decomp-gD}
 (i)
 $\hat\g^D = \gout^\trig \oplus \g^D_+ \oplus \C \hat k$.
 (ii)
 $M = V + \gout^\trig M$.
\end{lemma}

\begin{proof}
 (i)
 The space $\g^D$ is spanned by $\g^D_+$ and the elements of the form $0
 \oplus \cdots \oplus J_{a,b}\tensor\xi_i^{-n} \oplus \cdots \oplus 0$,
 where $(a,b) \in \{0,1,\dots,N-1\}^2 \setminus\{(0,0)\}$, $n\in
 \Z_{>0}$, $\xi_i$ is the local coordinate at $Q_i \in D$ and
 $J_{a,b}\tensor \xi_i^{-n}$ sits at the $i$-th component. The latter
 elements can be replaced with the following which belong to
 $\gout^\trig$: if $n=1$,
\begin{equation}
    J_{a,b}(t) :=
    \begin{cases}
    J_{a,b} t^a(t^N - t_i^N)^{-1} &(a\neq 0),
    \\
    J_{0,b} (\eps^b t^N - t_i^N)(t^N - t_i^N)^{-1} &(a=0).
    \end{cases}
\label{def:Jab(t)}
\end{equation}
 and if $n > 1$, $(t\, d/dt)^{n-1}J_{a,b}(t)$.  Thus we have
 $\gout^\trig + \g^D_+ = \g^D$. It is an exercise of complex analysis to
 show that $\gout^\trig \cap \g^D_+ = \{0\}$.

 (ii)
 By the assumption $M_i = U\hat\g^{Q_i} V_i$, $\g^{Q_i}_+ V_i \subset
 V_i$ and the Poincar\'e-Birkhoff-Witt theorem, we have $M = U\hat\g^D V
 = U\gout^\trig V = V + \gout^\trig M$.
\end{proof}
 
 By this lemma, every element $v\in M$ is decomposed as $v = v_f +
 v_\out$, $v_f \in V$, $v_\out \in \gout^\trig M$. The
 finite-dimensionality of $\CC_\trig(M)$ immediately follows from this
 expression. Note that $\gout^\trig$ acts trivially on $\univac_0$ and
 $\univac_\infty$. Hence $\univac_0 \tensor v_\out \tensor
 \univac_\infty$ belongs to $\gout^\orb(M(0)_0 \tensor M \tensor
 M(0)_\infty) \subset K$. $K$ becomes a right $U\h$-module by the action
 $v_0 \tensor v \tensor v_\infty \mapsto v_0 H \tensor v \tensor
 v_\infty$ of $H \in \h$. Therefore $p(H) \tensor v_\out \tensor
 \univac_\infty$ belongs to $K$, which allows us to assume that $v \in
 V$ in \eqref{p(H)*v*1<K}.

 We further reduces the problem and show \eqref{p(H)*v*1<K} for a weight
 vector in $V$, $v_\lambda \in V_\lambda$. $K$ being a right
 $U\h$-module, as we mentioned above, it is enough to show that
\begin{equation}
    (H_b - \tilde\lambda(H_b)) \tensor 
    v_\lambda \tensor \univac_\infty \in K,
\label{(H-lambda)*v*1<K}
\end{equation}
 for $b=1,\dots,N-1$. Since the constant function with value $H_b$
 belongs to $\gout^\orb$, we have
\begin{equation}
 \begin{split}
    H_b \tensor v_\lambda \tensor \univac_\infty
    &\equiv
    - \univac_0 \tensor H_b v_\lambda \tensor \univac_\infty
    - \univac_0 \tensor v_\lambda \tensor H_b
\\
    &=
    - \univac_0 \tensor \lambda(H_b) v_\lambda \tensor \univac_\infty
    - \univac_0 \tensor v_\lambda \tensor \Ad\beta(\Ad\beta^{-1}(H_b))
\\
    &\equiv
    - \lambda(H_b) \univac_0 \tensor v_\lambda \tensor \univac_\infty
    + \Ad\beta^{-1}(H_b) \tensor v_\lambda \tensor \univac_\infty
\\
    &=
    - \lambda(H_b) \univac_0 \tensor v_\lambda \tensor \univac_\infty
    + \eps^{-b} H_b \tensor v_\lambda \tensor \univac_\infty
 \end{split}
\label{H*v*1}
\end{equation}
 modulo $K$, which proves \eqref{(H-lambda)*v*1<K} and consequently
 \eqref{I*M*M<K}. Similarly we can prove that $M(0)_0 \tensor M \tensor
 I' \subset K$ and the inclusion \eqref{kernel-inclusion} is proved.

 Now we have the following commutative diagram with exact rows.
\begin{equation}
\minCDarrowwidth 3mm
\begin{CD}
     @. 0 @. 0 @.  @.
  \\
   @. @VVV @VVV  @.  @.
  \\
  0
  @>>> \tilde I
  @=   \tilde I
  @>>> 0
  @.
  \\
  @. @VVV @VVV @VVV @.
  \\
  0
  @>>> K
  @>>> M(0)_0 \tensor M \tensor M(0)_\infty
  @>>> \CC_\trig(M)
  @>>> 0
  \\
  @. @VVV @VVV @VVV @.
  \\
  0
  @>>> \bar K
  @>>> \bar M(0)_0 \tensor M \tensor \bar M(0)_\infty
  @>>> \bar M(0)_0 \tensor M \tensor \bar M(0)_\infty/\bar K
  @>>> 0
  \\
  @. @VVV @VVV @VVV @.
  \\
  @. 0 @. 0 @. 0 @.
\end{CD}
\label{cd:9lemma}
\end{equation}
 Here $\tilde I := I \tensor M \tensor M(0)_\infty + M(0)_0 \tensor M
 \tensor I'$ and $\bar K := K/\tilde I$. (The exactness of the middle row is
 due to \thmref{thm:general-factorisation}.) Applying the nine lemma in
 the homological algebra to this diagram (or chasing the diagram), we
 have
\begin{equation}
    \CC_\trig(M)
    \simeqq
    \bar M(0)_0 \tensor M \tensor \bar M(0)_\infty/\bar K.
\label{fact:finite:temp1}
\end{equation}
 It is easy to see that
\begin{equation}
    \bar K = 
    \gout^\orb (\bar M(0)_0 \tensor M \tensor \bar M(0)_\infty)
    + 
    (\bar M(0)_0 \tensor M \tensor \bar M(0)_\infty) \rho_{1,\beta}(\h),
\label{def:barK}
\end{equation}
 because the left $\hat\g^D$-action and the right $\h$-action on $M(0)_0
 \tensor M \tensor M(0)_\infty$ commute and the right $\h$-action is
 commutative. By \eqref{fact:finite:temp1} and \eqref{def:barK} we have
\begin{equation}
 \begin{split}
    \CC_\trig(M)
    &\simeqq
    \CC_\orb(\bar M(0)_0\tensor M \tensor \bar M(0)_\infty)/
    \overline{\CC_\orb}\cdot \rho_{1,\beta}(\h),
\\
    \overline{\CC_\orb}\cdot \rho_{1,\beta}(\h)
    &:=
    \CC_\orb(\bar M(0)_0\tensor M \tensor \bar M(0)_\infty) 
    \rho_{1,\beta}(\h).
 \end{split}
\label{fact:finite:temp2}
\end{equation}

\begin{lemma}
\label{lem:barM(0)}
\begin{equation}
    \bar M(0)_0 \simeqq 
    \bigoplus_{\lambda \in \wt(V)} M_{\tilde \lambda}, \qquad
    \bar M(0)_\infty \simeqq 
    \bigoplus_{\lambda \in \wt(V)} M_{\tilde \lambda'}.
\label{barM(0):decomp}
\end{equation}
(cf.\ \eqref{def:tilde-lambda}.)
\end{lemma}

\begin{proof}

 For two distinct weights $\lambda$ and $\mu$ in $\wt(V)$, we have
 $I_\lambda + I_\mu = U\h$, since $\{H_b\}_{b=1,\dots,N-1}$ is a basis
 of $\h$ and the map $\lambda \mapsto \tilde\lambda$ is injective. (See
 \eqref{def:Ih,I'h}.) Hence the ``Chinese remainder theorem'' implies
\begin{equation}
    U\h/I_\h \simeqq \bigoplus_{\lambda\in\wt(V)} U\h/I_\lambda.
\label{chinese-remainder-thm}
\end{equation}
 Regarding both hand sides as $\bigl(\g^{(0)}_+ \oplus \h \oplus \C \hat
 k\bigr)$-modules, on which $\g^{(0)}_+$ acts trivially and $\hat k$
 acts as $k$, we induce them up to $\hat\g^{(0)}$-modules and we obtain
 the first isomorphism in \eqref{barM(0):decomp}. The second isomorphism
 in \eqref{barM(0):decomp} is proved in the same manner.
\end{proof}

 Due to \lemref{lem:barM(0)} we may rewrite \eqref{fact:finite:temp2} as
\begin{equation}
 \begin{split}
    &\CC_\trig(M)
    \simeqq
    \bigoplus_{\lambda,\mu \in \wt(V)}
    \CC_\orb(M_{\tilde\lambda} \tensor M \tensor M_{\tilde\mu'})/
    \CC_{\orb,\lambda,\mu}\cdot \rho_{1,\beta}(\h),
\\
    &\CC_{\orb,\lambda,\mu}\cdot \rho_{1,\beta}(\h)
    :=
    \CC_\orb(M_{\tilde\lambda} \tensor M \tensor M_{\tilde\mu'})
    \rho_{1,\beta}(\h).
 \end{split}
\label{fact:finite:temp3}
\end{equation}
 We show that the component $\CC_\orb(M_{\tilde\lambda} \tensor M
 \tensor M_{\tilde\mu'})/ {\CC_{\orb,\lambda,\mu}}\cdot
 \rho_{1,\beta}(\h)$ of the above decomposition vanishes unless $\lambda
 = \mu$. Let us take $v_0 \tensor v \tensor v_\infty \in
 M_{\tilde\lambda} \tensor M \tensor M_{\tilde\mu'}$. For any $H\in\h$,
 we have $v_0 H = \tilde\lambda(H) v_0$ and $v_\infty H = \tilde\mu'(H)
 v_\infty$, which lead to
\begin{equation}
 \begin{split}
    (v_0 \tensor v \tensor v_\infty) \rho_{1,\beta}(H)
    &=
    v_0 H \tensor v \tensor v_\infty +
    v_0   \tensor v \tensor v_\infty \Ad(\beta) H
\\
    &=
    (\lambda - \mu) (1-\Ad\beta)^{-1} (H) 
    (v_0 \tensor v \tensor v_\infty),
 \end{split}
\label{(v0*v*voo)H}
\end{equation}
 due to \eqref{def:tilde-lambda}. Hence if $\lambda - \mu \neq 0$, there
 exists $H \in \h$ such that $(\lambda - \mu) (1-\Ad\beta)^{-1} (H) \neq
 0$, which means $v_0 \tensor v \tensor v_\infty \in
 \CC_{\orb,\lambda,\mu}\cdot \rho_{1,\beta}(\h)$. Therefore 
\begin{multline}
    \CC_\orb(M_{\tilde\lambda} \tensor M \tensor M_{\tilde\mu'})/
    \CC_{\orb,\lambda,\mu}\cdot \rho_{1,\beta}(\h)
\\
    =
    \begin{cases}
    \CC_\orb(M_{\tilde\lambda} \tensor M \tensor M_{\tilde\lambda'})
    &(\lambda = \mu),
   \\
    0 &(\lambda \neq \mu).
    \end{cases}
\label{fact:finite:component}
\end{multline}
 This completes the proof of \thmref{thm:factorisation:finite}.
\end{proof}

\section{Example (case of Weyl modules)}
\label{sec:weyl-mod}

Let us directly compute the space of conformal coinvariants for Weyl
modules as a simple example and see how the factorisation theorem
\thmref{thm:factorisation:finite} decompose $\CC_\trig$ into
$\CC_\orb$'s.

Let $V_i$ be a finite-dimensional $\g$-module and regard them as
$(\g^{Q_i}_+\oplus\C\hat k)$-module of level $k$ as before. The {\em
Weyl module} is the $\hat\g^{Q_i}$-module induced up from $V_i$:
\begin{equation}
    M(V_i):= \Ind_{\hat\g^{Q_i}_+}^{\g^{Q_i}\oplus\C\hat k} V_i.
\label{def:weyl}
\end{equation}
The tensor product $M(V) = M(V_1) \tensor \cdots \tensor M(V_L)$ is
naturally induced up from $V = V_1 \tensor\cdots \tensor V_L$ and, by
virtue of \lemref{lem:decomp-gD} (i), we have
\begin{equation}
    M(V) = \Ind_{\hat\g^D_+}^{\g^D\oplus\C\hat k} V
         = U\gout^\trig \tensor_\C V
         = V \oplus \gout^\trig\,  U\gout^\trig \tensor_\C V
         = V \oplus \gout^\trig\, M(V).
\label{M(V):decomp}
\end{equation}
Hence the space of conformal coinvariants is isomorphic to $V$ itself:
\begin{equation}
    \CC_\trig(M(V)) \simeqq V.
\label{CC(M(V))=V}
\end{equation}

On the other hand, let us compute the component
$\CC_\orb(M_{\tilde\lambda} \tensor M(V) \tensor M_{\tilde\lambda'})$ of
the factorisation, \eqref{factorisation:finite}. Put $\germa_1 =
\gout^\orb$, $\germa_2 = \borel^{(0)} \oplus \borel^D \oplus
\borel^{(\infty)} \oplus \C\hat k$, where $\borel^{(0)}$ and
$\borel^{(\infty)}$ are subalgebras of $\hat\g^{(0)}$ and
$\hat\g^{(\infty)}$ defined at the beginning of
\secref{sec:factorisation1} and $\borel^D := \oplus_{i=1}^L
\g^{Q_i}_+$. Then
\begin{equation}
     \germa := \germa_1 + \germa_2 = \hat\g^{D + (0) + (\infty)}, \qquad
     \germa_1 \cap \germa_2 = \h.
\label{a,a1*a2-for-weyl}
\end{equation}
We apply \lemref{lem:KL} to these algebras and the $\germa_2$-module $V
\simeqq \C 1_{\tilde\lambda} \tensor V \tensor \C
1_{\tilde\lambda'}$. The module $W$ in the lemma is
\begin{equation}
    \Ind_{\germa_2}^\germa 
    \C 1_{\tilde\lambda} \tensor V \tensor \C 1_{\tilde\lambda'}
    =
    M_{\tilde\lambda} \tensor M(V) \tensor M_{\tilde\lambda'}.
\label{W-for-weyl}
\end{equation}
The space of coinvariants ``$V/\germa_1 \cap \germa_2 V$'' in the lemma
is
\begin{equation}
    \C 1_{\tilde\lambda} \tensor V \tensor \C 1_{\tilde\lambda'}/
    \h (\C 1_{\tilde\lambda} \tensor V \tensor \C 1_{\tilde\lambda'}).
\label{V/aaV-for-weyl:temp}
\end{equation}
Take a weight vector $v_\mu \in V$ of weight $\mu$. The action of
$H\in\h$ on $1_{\tilde\lambda} \tensor v_\mu \tensor 1_{\tilde\lambda'}$
is equal to the multiplication of $\tilde\lambda(H) + \mu(H) +
\tilde\lambda'(H) = (-\lambda + \mu)(H)$. (See
\eqref{def:tilde-lambda}.) Namely, if $\lambda = \mu$,
$\h(1_{\tilde\lambda} \tensor v_\mu \tensor 1_{\tilde\lambda'}) = 0$,
while if $\lambda\neq \mu$, $\h (\C 1_{\tilde\lambda} \tensor V \tensor
\C 1_{\tilde\lambda'})$ contains $1_{\tilde\lambda} \tensor v_\mu
\tensor 1_{\tilde\lambda'}$. Consequently, 
\begin{equation}
    \C 1_{\tilde\lambda} \tensor V \tensor \C 1_{\tilde\lambda'}/
    \h (\C 1_{\tilde\lambda} \tensor V \tensor \C 1_{\tilde\lambda'})
    \simeqq
    V_\lambda.
\label{V/aaV-for-weyl}
\end{equation}
The coinvariants ``$W/\germa_1 W$'' in \lemref{lem:KL} is
\begin{equation}
    M_{\tilde\lambda} \tensor M(V) \tensor M_{\tilde\lambda'}/
    \gout^\orb M_{\tilde\lambda} \tensor M(V) \tensor M_{\tilde\lambda'}
    = \CC_\orb(M_{\tilde\lambda} \tensor M(V) \tensor M_{\tilde\lambda'}).
\label{W/a1W-for-weyl}
\end{equation}
because of \eqref{W-for-weyl}. Thus we have proved
\begin{equation}
    \CC_\orb(M_{\tilde\lambda} \tensor M(V) \tensor M_{\tilde\lambda'})
    \simeqq
    V_\lambda.
\label{fact:weyl:component}
\end{equation}
Summarising, the factorisation \eqref{factorisation:finite} for the Weyl
modules $M_i = M(V_i)$ is the same as the weight space decomposition of
$V$, $V = \bigoplus_{\lambda \in \wt(V)} V_\lambda$.

\section{Summary and concluding remarks}
\label{summary+remarks}

We formulated the twisted WZW model on a degenerate elliptic curve
$E_\trig$ and the orbifold $E_\orb$,
\defref{def:CC,CB:trig,orb}. \thmref{thm:general-factorisation} asserts
that the space of conformal coinvariants (and thus the space of
conformal blocks) on $E_\trig$ factorises into that on $E_\orb$. This
factorisation theorem holds for any smooth modules, but the result can
be refined for modules generated by $\g$-modules,
\thmref{thm:factorisation:finite}. The factorisation for the case of
Weyl modules is nothing more than the weight space decomposition of the
tensor product of $\g$-module, \eqref{CC(M(V))=V},
\eqref{fact:weyl:component}.

From these results arise several problems:
\begin{itemize}
 \item In \Refcite{tuy:89} the factorisation of the space of conformal
       coinvariants/blocks on singular curves together with the
       existence of a flat connection (the KZ connection) leads to the
       locally freeness of the sheaf of conformal
       coinvariants/blocks. One of our next tasks is to follow this line
       to analyse the structure of the KZ connection, which might help
       finding the integral representation of the solutions of the
       elliptic KZ equations.

 \item The weights $\tilde\lambda$ and $\tilde\lambda'$ in
       \thmref{thm:factorisation:finite} look quite strange. For
       example, $\tilde\lambda = - \lambda/2$ when $\g=sl_2(\C)$ ($N=2$)
       and $M_{\tilde\lambda}$ can hardly be integrable when $M_\lambda$
       is. But such weights are inevitable due to the charge
       conservation. The same kind of phenomenon is implicit but exists
       also in \Refcite{eti:94}, where the solution of the elliptic KZ
       equations is constructed as a twisted trace of vertex
       operators. Representation theoretical meaning of
       $M_{\tilde\lambda}$ is in question.

 \item We computed the space of conformal coinvariants for the Weyl
       modules in \secref{sec:weyl-mod}. This case is important in the
       study of the elliptic Gaudin model.\cite{kur-tak:97} The case
       for integrable modules instead of Weyl modules would also be
       important. 
\end{itemize}

\section*{Acknowledgments}

The author expresses his gratitude to Akihiro Tsuchiya who explained
details of his works, Michio Jimbo, Tetsuji Miwa, Hiroyuki Ochiai,
Kiyoshi Ohba, Nobuyoshi Takahashi and Tomohide Terasoma for discussion
and comments. He also thanks Gen Kuroki for discussions during the
previous collaborations. Christian Schweigert and Martin Halpern who
informed the references for the orbifold WZW model are gratefully
acknowledged.

The atmosphere and environment of Mathematical Scienences Research
Institute (Berkeley, U.S.A.), Institute for Theoretical and Experimental
Physics (Moscow, Russia) and the VIth International Workshop ``Conformal
Field Theory and Integrable Models'' (Chernogolovka, Russia) were very
important. The author thanks their hospitality.

This work is partly supported by the Grant-in-Aid for Scientific
Research, Japan Society for the Promotion of Science.

\end{document}